\documentclass{amsart}%
\usepackage{amssymb}
\usepackage{amsfonts}
\usepackage{geometry}
\usepackage{color}
\usepackage{amsmath}
\usepackage{verbatim}
\usepackage{graphicx}%
\providecommand{\U}[1]{\protect\rule{.1in}{.1in}}
\newtheorem{theorem}{Theorem}
\theoremstyle{plain}

\newtheorem{remark}[theorem]{Remark}

\numberwithin{equation}{section}
\begin{document}
\title{The $T1\,$theorem for the Hilbert transform fails when $p\not =2$}
\author[M. Alexis]{Michel Alexis}
\address{Department of Mathematics \& Statistics, McMaster University, 1280 Main Street
West, Hamilton, Ontario, Canada L8S 4K1}
\email{alexism@mcmaster.ca}
\author[J. L. Luna-Garcia]{Jose Luis Luna-Garcia}
\address{Department of Mathematics \& Statistics, McMaster University, 1280 Main Street
West, Hamilton, Ontario, Canada L8S 4K1}
\email{lunagaj@mcmaster.ca}
\author[E.T. Sawyer]{Eric T. Sawyer}
\address{Department of Mathematics \& Statistics, McMaster University, 1280 Main Street
West, Hamilton, Ontario, Canada L8S 4K1 }
\email{Sawyer@mcmaster.ca}
\author[I. Uriarte-Tuero]{Ignacio Uriarte-Tuero}
\address{Department of Mathematics, University of Toronto\\
Room 6290, 40 St. George Street, Toronto, Ontario, Canada M5S 2E4\\
(Adjunct appointment)\\
Department of Mathematics\\
619 Red Cedar Rd., room C212\\
Michigan State University\\
East Lansing, MI 48824 USA}
\email{ignacio.uriartetuero@utoronto.ca}
\thanks{E. Sawyer is partially supported by a grant from the National Research Council
of Canada}
\thanks{I. Uriarte-Tuero has been partially supported by grant MTM2015-65792-P
(MINECO, Spain), and is partially supported by a grant from the National
Research Council of Canada}

\begin{abstract}
Given $1<p<\infty$, $p\neq2$, we show that the $T1$ theorem for the Hilbert
transform fails for $L^{p}$, despite holding for $p=2$. More precisely, we
construct a pair of locally finite positive Borel measures $\left(
\sigma,\omega\right)  $ that satisfy the two-tailed $A_{p}$ condition, and
satisfy both of the $L$\thinspace$^{p}$-testing conditions for the Hilbert
transform $H$, yet $H_{\sigma}:L^{p}\left(  \sigma\right)  \not \rightarrow
L^{p}\left(  \omega\right)  $. In the opposite direction, the $T1$ theorem for
the Hilbert transform for $p=2$ was proved a decade ago in the two part paper
\cite{LaSaShUr3};\cite{Lac}.

\end{abstract}
\maketitle

\section{Introduction}

The $T1$ conjecture of F. Nazarov, S. Treil, and A. Volberg on the boundedness
of the Hilbert transform from one weighted space $L^{2}\left(  \sigma\right)
$ to another $L^{2}\left(  \omega\right)  $, was settled affirmatively in the
two part paper \cite{LaSaShUr3},\cite{Lac}, with the assumption of common
point masses subsequently removed by Hyt\"{o}nen in \cite{Hyt2}. One might
conjecture that the natural extension of this result to $p\neq2$ holds, namely
that the Hilbert transform $H$ is bounded from $L^{p}\left(  \sigma\right)  $
to $L^{p}\left(  \omega\right)  $ with general locally finite positive Borel
measures $\sigma$ and $\omega$ \emph{if and only if} the local testing
conditions for $H$,%
\begin{equation}
\int_{I}\left\vert H\mathbf{1}_{I}\sigma\right\vert ^{p}d\omega\leq
\mathfrak{T}_{H,p}\left(  \sigma,\omega\right)  ^{p}\left\vert I\right\vert
_{\sigma}\text{ and }\int_{I}\left\vert H\mathbf{1}_{I}\omega\right\vert
^{p^{\prime}}d\sigma\leq\mathfrak{T}_{H,p^{\prime}}\left(  \omega
,\sigma\right)  ^{p^{\prime}}\left\vert I\right\vert _{\omega},
\label{testing conditions}%
\end{equation}
both hold, along with the two-tailed Muckenhoupt $\mathcal{A}_{p}$ condition,%
\begin{equation}
\mathcal{A}_{p}\left(  \sigma,\omega\right)  \equiv\left(  \int_{\mathbb{R}%
}\frac{\left\vert I\right\vert ^{p-1}}{\left[  \left\vert I\right\vert
+\operatorname*{dist}\left(  x,I\right)  \right]  ^{p}}d\omega\right)
^{\frac{1}{p}}\left(  \int_{\mathbb{R}}\frac{\left\vert I\right\vert
^{p^{\prime}-1}}{\left[  \left\vert I\right\vert +\operatorname*{dist}\left(
x,I\right)  \right]  ^{p^{\prime}}}d\sigma\right)  ^{\frac{1}{p^{\prime}}%
}<\infty. \label{Muck conditions}%
\end{equation}

The purpose of this paper is to dash this hope by constructing a pair of
locally finite positive Borel measures $\left(  \omega,\sigma_{p}\right)  $ as
in \cite{LaSaUr2}, but with $\sigma_{p}$ adjusted according to $p$, for which
(\ref{testing conditions}) and (\ref{Muck conditions}) both hold, yet
$H_{\sigma}$ is unbounded, i.e. $\int_{\mathbb{R}}\left\vert Hf\sigma
\right\vert ^{p}d\omega=\infty$ and $\int_{\mathbb{R}}\left\vert f\right\vert
^{p}d\sigma<\infty$.

\begin{remark}
On the other hand, T. Hyt\"{o}nen and E. Vuorinen \cite{HyVu} have recently
conjectured an extension to weighted $L^{p}$ spaces that involves
\emph{quadratic} testing and Muckenhoupt conditions, obtained from testing the
$\ell^{2}$-valued extension of the norm inequality for $H$. In our
counterexample here, we derive the unboundedness of the Hilbert transform from
the failure of their quadratic Muckenhoupt condition. In the special case of
\emph{doubling measures}, their conjecture was recently proved for the
collection of vector Riesz transforms by the third author and B. Wick
\cite{SaWi}, with in fact the usual \emph{scalar} testing conditions in place
of the quadratic testing conditions - but with quadratic Muckenhoupt and weak
boundedness conditions similar to those in \cite{HyVu}.
\end{remark}

\begin{remark}
It was pointed out in a seminal paper \cite[page 911]{NTV} on the $T1$ theorem
in $L^{2}$ spaces for well-localized operators, that for \emph{general Haar
multipliers}, the two weight $T1$ theorem fails for all $p\neq2$, and a
Bellman proof can be found in an unpublished note of F. Nazarov \cite{Naz}.
However, Vuorinen has proved the `\emph{quadratic}' $T1$ theorem for dyadic
shifts in $L^{p}$ spaces \cite{Vuo}.
\end{remark}

Denote by $H$ the Hilbert transform of a complex Borel measure $\nu$ by%
\[
H\nu\left(  x\right)  \equiv\int_{\mathbb{R}}\frac{1}{y-x}d\nu\left(
y\right)  .
\]
Here is the \emph{quadratic} $\mathcal{A}_{p}$ condition of Hyt\"{o}nen and
Vuorinen \cite{HyVu},%
\begin{equation}
\left\Vert \left(  \sum_{i=1}^{\infty}\left(  \mathbf{1}_{I_{i}}%
\int_{\mathbb{R}\setminus I_{i}}\frac{f_{i}\left(  x\right)  }{\left\vert
x-c_{i}\right\vert }d\sigma\left(  x\right)  \right)  ^{2}\right)  ^{\frac
{1}{2}}\right\Vert _{L^{p}\left(  \omega\right)  }\leq\mathcal{A}_{p}%
^{\ell^{2}}\left(  \sigma,\omega\right)  \left\Vert \left(  \sum_{i=1}%
^{\infty}f_{i}^{2}\right)  ^{\frac{1}{2}}\right\Vert _{L^{p}\left(
\sigma\right)  }, \label{quad A2 cond}%
\end{equation}
taken over all sequences of intervals $\left\{  I_{i}\right\}  _{i=1}^{M}$,
and all sequences of functions $\left\{  f_{i}\right\}  _{i=1}^{M}$, and all
$M\in\mathbb{N}$. The constant $\mathcal{A}_{p}^{\ell^{2}}\left(
\sigma,\omega\right)  $ is by definition the best possible. As pointed out in
\cite{HyVu}, this condition is necessary for the norm inequality. We will
construct our counterexample $\left(  \sigma,\omega\right)  $ to the $T1$
theorem for $H$ so as to satisfy the $A_{p}$ and testing conditions for
$\left(  \sigma,\omega\right)  $, but \emph{fail} the quadratic $\mathcal{A}%
_{p}$ condition.

\section{Cantor construction}

We will use analogues of the measure pairs $\left(  \sigma,\omega\right)  $
and $\left(  \dot{\sigma},\omega\right)  $ constructed in \cite{LaSaUr2} (that
were used to demonstrate the failure of necessity of the pivotal condition).
Recall the middle-third Cantor set $\mathsf{E}$ and Cantor measure $\omega$ on
the closed unit interval $I_{1}^{0}=\left[  0,1\right]  $. At the $k^{th}$
generation in the construction, there is a collection $\left\{  I_{j}%
^{k}\right\}  _{j=1}^{2^{k}}$ of $2^{k}$ pairwise disjoint closed intervals of
length $\left\vert I_{j}^{k}\right\vert =\frac{1}{3^{k}}$. With $K_{k}%
=\bigcup_{j=1}^{2^{k}}I_{j}^{k}$, the Cantor set is defined by $\mathsf{E}%
=\bigcap_{k=1}^{\infty}K_{k}=\bigcap_{k=1}^{\infty}\left(  \bigcup
_{j=1}^{2^{k}}I_{j}^{k}\right)  $. The Cantor measure $\omega$ is the unique
probability measure supported in $\mathsf{E}$ with the property that it is
equidistributed among the intervals $\left\{  I_{j}^{k}\right\}  _{j=1}%
^{2^{k}}$ at each scale $k$, i.e.%
\begin{equation}
\left\vert I_{j}^{k}\right\vert _{\omega}=2^{-k},\ \ \ \ \ k\geq0,1\leq
j\leq2^{k}. \label{omega measure}%
\end{equation}

Denote the removed open middle third of $I_{j}^{k}$ by $G_{j}^{k}$. Let
$\dot{z}_{j}^{k}\in G_{j}^{k}$ be the center of the interval $G_{j}^{k}$,
which is also the center of the interval $I_{j}^{k}$, and for $1<p<\infty$,
define
\begin{equation}
\label{def_sigma_dot}\dot{\sigma}_{p}=\sum_{k,j}s_{j}^{k}\delta_{\dot{z}%
_{j}^{k}},
\end{equation}
where the sequence of positive numbers $s_{j}^{k}$ is chosen to satisfy the
following precursor of the $A_{p}$ condition:%
\[
\frac{\left(  s_{j}^{k}\right)  ^{p-1}\left\vert I_{j}^{k}\right\vert
_{\omega}}{\left\vert I_{j}^{k}\right\vert ^{p}}=1,\ \ \ s_{j}^{k}%
=2^{-k}\left(  \frac{2}{3}\right)  ^{kp^{\prime}}\ \ \ \ \ k\geq0,1\leq
j\leq2^{k}.
\]

To define the measure $\sigma_{p}$ for the counter example, we modify
$\dot{\sigma}_{p}$ as follows, and for convenience in notation we will drop
the subscript $p$ from both $\dot{\sigma}_{p}$ and $\sigma_{p}$. Shift the
point masses $s_{j}^{k}\delta_{\dot{z}_{j}^{k}}$ in $\dot{\sigma}$ to new
positions $s_{j}^{k}\delta_{z_{j}^{k}}$ in $G_{j}^{k}$ where $z_{j}^{k}$ is
the unique point $z$ in the open interval $G_{j}^{k}$ where $H\omega\left(
z\right)  $ vanishes. Such a point $z_{j}^{k}$ exists because $H\omega\left(
z\right)  $ is strictly increasing from $-\infty$ to $+\infty$ on the interval
$G_{j}^{k}$. Let
\begin{equation}
\label{def_sigma}\sigma=\sum_{k,j}s_{j}^{k}\delta_{z_{j}^{k}}.
\end{equation}

In the case $p=2$, it was shown in \cite{LaSaUr2} that the two-tailed
Muckenhoupt condition (\ref{Muck conditions}), and the forward and backward
testing conditions (\ref{testing conditions}) hold for both of the measure
pairs $\left(  \dot{\sigma},\omega\right)  $ and $\left(  \sigma
,\omega\right)  $. These proofs extend easily to the case $p\neq2$, and are
left for the reader. For example, the self-similarity argument in
\cite{LaSaUr2} gives
\begin{equation}
\int\left\vert H\dot{\sigma}\right\vert ^{p}\omega=\frac{2^{p^{\prime}-1}%
}{3^{p^{\prime}}}\left(  1+\varepsilon\right)  \int\left\vert H\dot{\sigma
}\right\vert ^{p}\omega+\frac{2^{p^{\prime}-1}}{3^{p^{\prime}}}\left(
1+\varepsilon\right)  \int\left\vert H\dot{\sigma}\right\vert ^{p}%
\omega+\mathcal{R}_{\varepsilon}, \label{reproduce}%
\end{equation}
where $\mathcal{R}_{\varepsilon}\leq C_{\varepsilon}\mathcal{A}_{p}\left(
\int\dot{\sigma}\right)  $, and hence%
\[
\int\left\vert H\dot{\sigma}\right\vert ^{p}\omega=\frac{1}{1-2\frac
{2^{p\prime-1}}{3^{p^{\prime}}}\left(  1+\varepsilon\right)  }\mathcal{R}%
_{\varepsilon}\leq C_{\varepsilon}\mathcal{A}_{p}^{2}\left(  \int\dot{\sigma
}\right)  ,
\]
provided $0<\varepsilon<\left(  \frac{3}{2}\right)  ^{p^{\prime}}-1$.

\subsection{Failure of $\mathcal{A}_{p}^{\ell^{2}}\left(  \sigma
,\omega\right)  $}

It remains to show that $\mathcal{A}_{p}^{\ell^{2}}\left(  \sigma
,\omega\right)  =\infty$ for $p\neq2$. We will show that there is $f=\left\{
f_{i}\right\}  _{i=1}^{\infty}$ satisfying%
\begin{equation}
\left\Vert \left(  \sum_{i=1}^{\infty}f_{i}^{2}\right)  ^{\frac{1}{2}%
}\right\Vert _{L^{p}\left(  \sigma\right)  }^{p}<\infty\text{ and }\left\Vert
\left(  \sum_{i=1}^{\infty}\left(  \mathbf{1}_{I_{i}}\int_{\mathbb{R}%
\setminus3I_{i}}\frac{f_{i}\left(  x\right)  }{\left\vert x-c_{i}\right\vert
}d\sigma\left(  x\right)  \right)  ^{2}\right)  ^{\frac{1}{2}}\right\Vert
_{L^{p}\left(  \omega\right)  }^{p}=\infty. \label{sup rat}%
\end{equation}
In fact we will take
\[
\left\{  I_{i}\right\}  _{i=1}^{\infty}\equiv\left\{  I_{j}^{k}\right\}
_{k,j}\text{ and }f_{i}=f_{k,j}\equiv b_{j}^{k}\mathbf{1}_{\theta I_{j}^{k}%
}\ ,
\]
where $\theta I_{j}^{k}$ denotes the triadic sibling of $I_{j}^{k}$, and use
the estimates
\begin{align*}
&  \left\vert I_{j}^{k}\right\vert =3^{-k},\ \ \ \left\vert I_{j}%
^{k}\right\vert _{\omega}=2^{-k},\ \ \ ,\left\vert I_{j}^{k}\right\vert
_{\sigma}\approx\widehat{s}_{j}^{k}\left[  p\right]  =\frac{2^{k\left(
p^{\prime}-1\right)  }}{3^{kp^{\prime}}},\\
&  \int_{\mathbb{R}\setminus3I_{j}^{k}}\frac{f_{k,j}\left(  x\right)
}{\left\vert x-c_{i}\right\vert }d\sigma\left(  x\right)  \approx b_{j}%
^{k}\int_{\theta I_{j}^{k}}\frac{1}{\left\vert I_{j}^{k}\right\vert }%
d\sigma\left(  x\right)  \approx b_{j}^{k}\frac{\left\vert I_{j}%
^{k}\right\vert _{\sigma}}{\left\vert I_{j}^{k}\right\vert }\approx b_{j}%
^{k}\frac{2^{k\left(  p^{\prime}-1\right)  }}{3^{k\left(  p^{\prime}-1\right)
}}.
\end{align*}

We define%
\[
b_{j}^{k}=\beta_{k}=\frac{1}{\left(  k+1\right)  ^{\frac{1}{p}}\left(
\ln\left(  k+2\right)  \right)  ^{\frac{1+\delta}{p}}}\left(  \frac{3}%
{2}\right)  ^{\left(  p^{\prime}-1\right)  k}\equiv\frac{1}{a_{k,p,\delta}%
}\left(  \frac{3}{2}\right)  ^{\left(  p^{\prime}-1\right)  k},
\]
and obtain%
\[
\left\Vert \left(  \sum_{i=1}^{\infty}f_{i}^{2}\right)  ^{\frac{1}{2}%
}\right\Vert _{L^{p}\left(  \sigma\right)  }^{p}=\int_{\mathbb{R}}\left(
\sum_{k,j}^{\infty}\left(  \beta_{k}\mathbf{1}_{\theta I_{j}^{k}}\left(
y\right)  \right)  ^{2}\right)  ^{\frac{p}{2}}d\sigma\left(  y\right)  .
\]
Set $\Omega_{k}\equiv\bigcup_{j=1}^{2^{k}}I_{j}^{k}=\bigcup_{j=1}^{2^{k}%
}\theta I_{j}^{k}$. We have the estimate%
\begin{align*}
\sum_{k,j}^{\infty}\beta_{k}^{2}\mathbf{1}_{\theta I_{j}^{k}}\left(  y\right)
&  \approx\sum_{k=0}^{N}\beta_{k}^{2}\mathbf{1}_{\Omega_{k}}\left(  y\right)
=\sum_{k=0}^{N}\beta_{k}^{2}\sum_{\ell=k}^{\infty}\mathbf{1}_{\Omega_{\ell
}\setminus\Omega_{\ell+1}}\left(  y\right) \\
&  =\sum_{\ell=0}^{\infty}\left(  \sum_{k=1}^{\ell}\beta_{k}^{2}\right)
\mathbf{1}_{\Omega_{\ell}\setminus\Omega_{\ell+1}}\left(  y\right)
\approx\sum_{\ell=0}^{\infty}\beta_{\ell}^{2}\mathbf{1}_{\Omega_{\ell
}\setminus\Omega_{\ell+1}}\left(  y\right)  ,
\end{align*}
and so changing dummy variables,%
\begin{align*}
&  \int_{\mathbb{R}}\left(  \sum_{k,j}^{\infty}\left(  b_{j}^{k}%
\mathbf{1}_{\theta I_{j}^{k}}\left(  y\right)  \right)  ^{2}\right)
^{\frac{p}{2}}d\sigma\left(  y\right)  \approx\int_{\mathbb{R}}\left(
\sum_{k=0}^{\infty}\beta_{k}^{2}\mathbf{1}_{\Omega_{k}\setminus\Omega_{k+1}%
}\left(  y\right)  \right)  ^{\frac{p}{2}}d\sigma\left(  y\right) \\
&  \approx\sum_{k=0}^{\infty}\beta_{k}^{p}\left\vert \Omega_{k}\setminus
\Omega_{k+1}\right\vert _{\sigma}\approx\sum_{k=0}^{\infty}\beta_{k}^{p}%
\sum_{j=1}^{2^{k}}\left\vert I_{j}^{k}\setminus\Omega_{k+1}\right\vert
_{\sigma}=\sum_{k=0}^{\infty}\beta_{k}^{p}\sum_{j=1}^{2^{k}}s_{j}^{k}\\
&  =\sum_{k=0}^{\infty}\beta_{k}^{p}\sum_{j=1}^{2^{k}}\frac{2^{k\left(
p^{\prime}-1\right)  }}{3^{kp^{\prime}}}=\sum_{k=0}^{\infty}\left(  \frac
{1}{a_{k,p,\delta}}\left(  \frac{3}{2}\right)  ^{\left(  p^{\prime}-1\right)
k}\right)  ^{p}\frac{2^{kp^{\prime}}}{3^{kp^{\prime}}}=\sum_{k=0}^{\infty
}\frac{1}{\left(  k+1\right)  \left(  \ln\left(  k+2\right)  \right)
^{1+\delta}}<\infty.
\end{align*}

Next we compute%
\begin{align*}
&  \left\Vert \left(  \sum_{i}\left(  \mathbf{1}_{I_{i}}\int_{\mathbb{R}%
\setminus3I_{i}}\frac{f_{i}\left(  x\right)  }{\left\vert x-c_{i}\right\vert
}d\sigma\left(  x\right)  \right)  ^{2}\right)  ^{\frac{1}{2}}\right\Vert
_{L^{p}\left(  \omega\right)  }^{p}\approx\int_{\mathbb{R}}\left(  \sum
_{k,j}\left(  \mathbf{1}_{I_{j}^{k}}\left(  y\right)  \beta_{k}\frac
{2^{k\left(  p^{\prime}-1\right)  }}{3^{k\left(  p^{\prime}-1\right)  }%
}\right)  ^{2}\right)  ^{\frac{p}{2}}d\omega\left(  y\right) \\
&  \approx\int_{\mathbb{R}}\left(  \sum_{k,j}\left(  \mathbf{1}_{I_{j}^{k}%
}\left(  y\right)  \frac{1}{a_{k,p,\delta}}\left(  \frac{3}{2}\right)
^{\left(  p^{\prime}-1\right)  k}\frac{2^{k\left(  p^{\prime}-1\right)  }%
}{3^{k\left(  p^{\prime}-1\right)  }}\right)  ^{2}\right)  ^{\frac{p}{2}%
}d\omega\left(  y\right)  =\int_{\mathbb{R}}\left(  \sum_{k,j}\left(
\mathbf{1}_{I_{j}^{k}}\left(  y\right)  \frac{1}{a_{k,p,\delta}}\right)
^{2}\right)  ^{\frac{p}{2}}d\omega\left(  y\right)  ,
\end{align*}
and truncating the sum in $k$ at $N$ we get%
\[
\sum_{k\leq N,j}\left(  \mathbf{1}_{I_{j}^{k}}\left(  y\right)  \frac
{1}{a_{k,p,\delta}}\right)  ^{2}=\sum_{k=0}^{N}\frac{1}{\left(  k+1\right)
^{\frac{2}{p}}\left(  \ln\left(  k+2\right)  \right)  ^{\frac{2}{p}\left(
1+\delta\right)  }}\mathbf{1}_{\Omega_{N}}\left(  y\right)  \approx
\frac{N^{1-\frac{2}{p}}}{\left(  \ln N\right)  ^{\frac{2}{p}\left(
1+\delta\right)  }}\mathbf{1}_{\Omega_{N}}\left(  y\right)  \text{ for }%
\omega\text{-a.e. }y,
\]
and so
\[
\int_{\mathbb{R}}\left(  \sum_{k\leq N,j}\left(  \mathbf{1}_{I_{j}^{k}}\left(
y\right)  \frac{1}{a_{k,p,\delta}}\right)  ^{2}\right)  ^{\frac{p}{2}}%
d\omega\left(  y\right)  \approx\int_{\mathbb{R}}\left(  \frac{N^{1-\frac
{2}{p}}}{\left(  \ln N\right)  ^{\frac{2}{p}\left(  1+\delta\right)  }%
}\mathbf{1}_{\Omega_{N}}\left(  y\right)  \right)  ^{\frac{p}{2}}%
d\omega\left(  y\right)  \approx\frac{N^{\frac{p}{2}-1}}{\left(  \ln N\right)
^{\left(  1+\delta\right)  }}\left\vert \Omega_{N}\right\vert _{\omega},
\]
which blows up to $\infty$ as $N\nearrow\infty$ provided $p>2$.

This shows that the $T1$ theorem fails for the Hilbert transform when $p>2$,
and a duality argument shows the same when $1<p<2$.

\begin{remark}
The arguments outlined above for the Hilbert transform generalize to
individual Riesz transforms $R_{j}$ in higher-dimensions. Indeed, if
considering $R_{1}$ in $\mathbb{R}^{n}$, take $\omega$ to the Cantor measure
on the segment $[0,1]\times\{0\}^{n-1}$, and take $\dot{\sigma}$ and $\sigma$
as in (\ref{def_sigma_dot}) and (\ref{def_sigma}) with $s_{j}^{k}\equiv
\frac{2^{k\left(  p^{\prime}-1\right)  }}{3^{nkp^{\prime}}}$, $\dot{z}_{j}%
^{k}$ as the centers of the middle-third line segments $G_{j}^{k}$ removed in
the construction of $\omega$, and $z_{j}^{k}$ defined as any one of the zeros
of $R_{1}\omega(x)$ occurring in $G_{j}^{k}$. Finally, to show the failure of
$\mathcal{A}_{p}^{\ell^{2}}(\sigma,\omega)$, follow the same proof as outlined
above, and take $b_{j}^{k}=\beta_{k}=\frac{1}{k^{\frac{1}{p}}\left(  \ln
k\right)  ^{\frac{1+\delta}{p}}\left(  \frac{3^{n}}{2}\right)  ^{(p^{\prime
}-1)k}}$. This implies the norm inequality fails for $R_{1}$: indeed, since
$\sigma$ and $\omega$ are supported on the line segment $[0,1]\times
\{0\}^{n-1}$, then \emph{ for this pair of measures} we can show that
$A_{p}^{\ell^{2}}$ is also necessary for the norm inequality for the Riesz
transform $R_{1}$.
\end{remark}


\begin{thebibliography}{999999999}                                                                                        %


\bibitem[Hyt2]{Hyt2}\textsc{Hyt\"{o}nen, Tuomas, }\textit{The two weight
inequality for the Hilbert transform with general measures,
\texttt{arXiv:1312.0843v2}.}

\bibitem[HyVu]{HyVu}\textsc{Tuomas Hyt\"{o}nen and Emil Vuorinen, }\textit{A
two weight inequality between }$L^{p}\left(  \ell^{2}\right)  $\textit{\ and
}$L^{p}$, \texttt{arXiv:1608.07199v2}.

\bibitem[Lac]{Lac}\textsc{Michael T. Lacey, }\textit{\ Two weight inequality
for the Hilbert transform: A real variable characterization, II}, Duke Math.
J. Volume \textbf{163}, Number 15 (2014), 2821-2840.

\bibitem[LaSaShUr3]{LaSaShUr3}\textsc{Michael T. Lacey, Eric T. Sawyer,
Chun-Yen Shen, and Ignacio Uriarte-Tuero,} \textit{Two weight inequality for
the Hilbert transform: A real variable characterization I}, Duke Math. J,
Volume \textbf{163}, Number 15 (2014), 2795-2820.

\bibitem[LaSaUr2]{LaSaUr2}\textsc{Lacey, Michael T., Sawyer, Eric T.,
Uriarte-Tuero, Ignacio,} \textit{A Two Weight Inequality for the Hilbert
transform assuming an energy hypothesis, } Journal of Functional Analysis,
Volume \textbf{263} (2012), Issue 2, 305-363.

\bibitem[Naz]{Naz}\textsc{F. Nazarov,} unpublished manuscript.

\bibitem[NTV]{NTV}\textsc{F. Nazarov, S. Treil and A. Volberg,} \textit{The
Bellman Functions and Two-Weight Inequalities for Haar Multipliers, }Journal
of the American Mathematical Society , Oct., 1999, Vol. \textbf{12}, No. 4
(Oct., 1999), 909-928.

\bibitem[SaWi]{SaWi}\textsc{Eric T. Sawyer and Brett D. Wick,} \textit{Two
weight }$L^{p}$\textit{ inequalities for }$\lambda$\textit{-fractional vector
Riesz tranforms and doubling measures}, \texttt{arXiv:2211.01920}.

\bibitem[Vuo]{Vuo}\textsc{Emil Vuorinen, }Two weight $L^{p}$%
-inequalities\textit{\ for dyadic shifts and the dyadic square function},
Studia Math. \textbf{237} (1) \ (2017), 25-56.
\end{thebibliography}
\end{document}